# Rectifiability of flat chains

By Brian White*

*In memory of Frederick J. Almgren Jr., 1933–1997*

## Abstract

We prove (without using Federer's structure theorem) that a finite-mass flat chain over any coefficient group is rectifiable if and only if almost all of its 0-dimensional slices are rectifiable. This implies that every flat chain of finite mass and finite size is rectifiable. It also leads to a simple necessary and sufficient condition on the coefficient group in order for every finite-mass flat chain to be rectifiable.

Let **G** be an abelian group with a norm that makes **G** a complete metric space. In 1966, Fleming [FL] developed a theory of flat chains with coefficients in **G**. In case **G** is the integers or the real numbers, with the usual norm, these flat chains are the integral or real flat chains, respectively, previously introduced by Federer and Fleming [FF]. In case **G** is the integers modulo $p$, these are the flat chains modulo $p$ which were later constructed in a different way by Federer [F]. After Fleming's paper, most research involving flat chains has been limited to integral, real, and modulo $p$ chains. However, recently it has become apparent that more general groups are of considerable interest. For example, other normed groups arise naturally in modelling immiscible fluids and soap bubble clusters [W2], and in proving that various surfaces are area-minimizing [ML].

Perhaps the most profound fact about integral flat chains is the rectifiability theorem: every integral flat chain of finite mass is rectifiable. (The compactness theorem for integral currents is an immediate consequence via

*The author was partially funded by NSF Grant DMS-95-04456.
1991 *Mathematics Subject Classification*. Primary 49Q15; secondary 49Q20.

[FL, 7.5].) The rectifiability theorem also holds for flat chains with coefficients in any finite group, but not for real flat chains. That raises the question: for which normed groups is every finite-mass flat chain rectifiable? In this paper we give a simple answer: those groups **G** that contain no nonconstant continuous path of finite length. Of course this includes any group (such as **Z** or any finite group) in which the norm of nonzero elements is bounded away from 0. But it also includes other interesting normed groups, such as:

(1) The $p$-adic numbers (where $p$ is a given prime), that is, the metric space completion of the rational numbers **Q** with respect to the metric given by the norm
$$\left| p^k \frac{m}{n} \right| = p^{-k}$$
where $p$ does not divide $m$ or $n$.

(2) The $p$-adic integers, that is the completion of **Z** under the same norm.

(3) The real numbers with the norm
$$|x| = \|x\|^\alpha$$
where $\|x\|$ is the usual absolute value of $x$ and $\alpha$ is a fixed number with $0 < \alpha < 1$.

Of course, **R** with its usual norm is not such a group. Consequently there exist finite mass real flat chains that are not rectifiable. However, we show that (for any coefficient group) every flat chain of finite mass and finite size is rectifiable. (Roughly speaking, mass and size are weighted area and unweighted area, respectively. See Section 8 for a precise definition.) In particular, since size is lower-semicontinuous, this implies Almgren's compactness theorem [A, 4.3] for real rectifiable currents of bounded size and mass.

The proofs in this paper are perhaps simpler than existing proofs of the rectifiability of finite mass integral flat chains. In particular, this paper does not use Federer's structure theorem for sets of finite Hausdorff measure. (Solomon [So], Almgren [A], and I [W1] have given other proofs that avoid the structure theorem, but [So] and [W1] use differential forms in an essential way and thus are limited to integral flat chains, and Almgren's proof requires lengthy preliminaries about multivalued functions.) Indeed, although the structure theorem was essential in the original proofs of rectifiability, I do not see how it can be used in proving the results of this paper.

The key theorem of this paper, from which the results mentioned above are derived, is the following:

RECTIFIABLE SLICES THEOREM. *Let A be a flat k-chain of finite mass in* $\mathbf{R}^N$. *The following are equivalent*:

(1) *A is rectifiable.*
(2) *For almost every $(N-k)$-plane $P$ parallel to a coordinate plane, the slice $A \cap P$ is a rectifiable $0$-chain.*

This theorem reduces questions about rectifiability of $k$-chains to much simpler questions about rectifiability of $0$-chains. Rectifiable $0$-chains are simply chains of the form $\sum_{i=1}^{\infty} g_i[x_i]$ where $x_i \in \mathbf{R}^N$ and $g_i \in \mathbf{G}$ with $\sum |g_i| < \infty$.

(In the case $\mathbf{G} = \mathbf{R}$, the rectifiable slices theorem has been discovered independently by L. Ambrosio and B. Kirchheim [AK]. Indeed, they develop a theory of rectifiable currents in metric spaces, and they prove an analog of the rectifiable slices theorem in that setting.)

The organization of this paper is as follows. In Section 1 we introduce notation and recall basic facts about flat chains and about approximate derivatives. In Section 2, we study 0-dimensional flat chains. In particular, we prove that there is a natural homomorphism $\chi$ from the group $\mathcal{F}_0(\mathbf{R}^N; \mathbf{G})$ of flat 0-chains to the coefficient group $\mathbf{G}$. The homomorphism $\chi$ acts on rectifiable 0-chains by summing the multiplicities: $\chi(\sum g_i[x_i]) = \sum g_i$. We also show that flat 0-chains of finite mass are equivalent to compactly supported $\mathbf{G}$-valued Borel measures of finite total variation.

In Section 3, we define slicing by affine planes and prove its basic properties. In particular, we prove (i) if $A$ is rectifiable, then almost all of its 0-dimensional slices are rectifiable, and (ii) if $A$ is a flat chain almost all of whose 0-dimensional slices vanish, then $A$ is 0. The first assertion is straightforward. The second is a direct consequence of the deformation theorem of [W3]. The idea of the proof of (ii) is the following. According to the deformation theorem, any flat chain $A$ can be approximated by a polyhedral chain $P$ in a cubical grid. This polyhedral chain depends only on the intersections (or slices) of $A$ with the various $(N-k)$-planes in the dual grid. The hypothesis in (ii) implies that if the grid is in general position, then $P$ will be 0. Since we get arbitrarily good approximations $P$ in this way, $A$ must be 0.

In Section 4, we use (ii) to reduce the rectifiable slices theorem to the special case of a $k$-chain that is supported in the graph of a continuous function $\Phi$ over a compact subset of a $k$-plane. In Section 5, we prove that, in case $k = 1$, the function $\Phi$ is approximately differentiable almost everywhere. In Section 6, for arbitrary $k$, we prove that the approximate partial derivative $\text{ap}\mathbf{D}_i \Phi$ exist almost everywhere; this is done by slicing $A$ to get a 1-dimensional chain and applying Section 4. It follows that $\Phi$ is approximately differentiable almost everywhere and hence that $A$ is rectifiable. That completes the proof of the rectifiable slices theorem.

In Section 7, we characterize those normed groups for which every finite mass flat chain is rectifiable.

Finally, in Section 8, we discuss size and other weighted area functionals.

# 1. Notation

Everything we use about flat chains is proved in the first half (through §6) of Fleming's paper [FL], except for the improved deformation theorem [W3]. Here we summarize the main facts about flat chains.

1.1. Let **G** be an abelian group. Let $|\cdot| : \mathbf{G} \to \mathbf{R}$ be a norm on **G**, that is, a function such that

(1) $|-g| \equiv |g|$,

(2) $|g + h| \leq |g| + |h|$, and

(3) $|g| \geq 0$ with equality if and only if $g = 0$.

We assume **G** is a complete metric space with respect to the induced metric $\operatorname{dist}(g, h) = |g - h|$. (Otherwise replace **G** by its metric space completion.)

Let $K$ be a closed convex subset of $\mathbf{R}^N$. We let $\mathcal{F}_k(K; \mathbf{G})$, $\mathcal{M}_k(K; \mathbf{G})$, and $\mathcal{P}_k(K; \mathbf{G})$ denote the spaces of flat $k$-chains, finite-mass flat $k$-chains, and polyhedral $k$-chains, respectively, in $K$ with coefficients in **G**. The spaces $\mathcal{F}_k(K; \mathbf{G})$ and $\mathcal{M}_k(K; \mathbf{G})$ are normed groups under the flat norm $A \mapsto \mathcal{F}(A)$ and the mass norm $A \mapsto \mathrm{M}(A)$, respectively; if $K$ is compact, then they are complete metric spaces. The flat norm and mass norms are related by

$$\mathcal{F}(A) = \inf \{\mathrm{M}(A - \partial B) + \mathrm{M}(B) : B \in \mathcal{F}_{k+1}\}.$$

If $A$ is polyhedral, then we may take the infimum over polyhedral $(k+1)$-chains $B$.

(If $K$ is compact, then $\mathcal{F}_k(K; \mathbf{G})$ is actually defined to be the $\mathcal{F}$-completion of $\mathcal{P}_k(K; \mathbf{G})$. Fleming also uses this as the definition when $K = \mathbf{R}^N$. However, I prefer to let $\mathcal{F}_k(\mathbf{R}^N; \mathbf{G})$ be the union of $\mathcal{F}_k(K; \mathbf{G})$ over compact sets $K$; this guarantees that every flat chain has compact support. The results of this paper are true with either definition, although the statement of subsection 2.2 needs to be modified if one allows noncompactly supported flat chains.)

Mass is lower-semicontinuous on $\mathcal{F}_k(K; \mathbf{G})$. Also, given any flat chain $A$, there exist polyhedral chains $A_i$ converging to $A$ with $\mathrm{M}(A) = \lim \mathrm{M}(A_i)$. A polyhedral chain $A$ can be expressed as a **G**-linear combination of a finite collection of non-overlapping oriented $k$-dimensional polyhedra:

$$A = \sum g_i [\sigma_i].$$

The mass of $A$ is then $\sum |g_i| \operatorname{area}(\sigma_i)$.

Associated to every finite mass $k$-chain $A$ and Borel set $S$, there is a finite mass $k$-chain $A \llcorner S$ that is, roughly speaking, the portion of $A$ in $S$. Then $S \mapsto A \llcorner S$ defines a $\mathcal{M}_k$-valued Borel measure, and $S \mapsto \mathrm{M}(A \llcorner S)$ defines a Radon measure $\mu_A$. (In Fleming's paper, $A \llcorner S$ is written $A \cap S$, but the notation we use has become standard. In this paper $A \cap (\cdot)$ denotes slicing; see §3.)

A Lipschitz map $f : \mathbf{R}^N \to \mathbf{R}^m$ induces a homomorphism $f_\# : \mathcal{F}_k(\mathbf{R}^N; \mathbf{G}) \to \mathcal{F}_k(\mathbf{R}^m; \mathbf{G})$. If $f$ is affine and $A$ is polyhedral, it is given by

$$f_\# : \sum g_i[\sigma_i] \mapsto \sum g_i[f(\sigma_i)].$$

The maps that occur in this paper are translations $\tau_y : x \mapsto x+y$, dilations $r : x \mapsto rx$, and projections onto subspaces. For aesthetic reasons, $\tau_{y\#}A$ will be written as $\tau_y A$.

1.2. *Rectifiable flat chains.* A set $S$ in $\mathbf{R}^N$ is said to be countably $k$-rectifiable if $S$ is contained in the countable union of Lipshitz images of $\mathbf{R}^k$. A finite-mass flat $k$-chain $A$ in $\mathbf{R}^N$ is said to be rectifiable if there exists a countably $k$-rectifiable Borel set $S$ such that $A = A \llcorner S$ (or, equivalently, such that $\mu_A(\mathbf{R}^N \setminus S) = 0$). Since $\mu_A V = 0$ for $\mathcal{H}^k$-nullsets $V$ (by [W3, §3]) we can choose the set $S$ to be a countable union of $k$-dimensional $C^1$ submanifolds of $\mathbf{R}^N$ (cf. [S, 11.1].)

Consequently, it is not hard to show that $A$ is rectifiable if and only if for every $\epsilon > 0$ there is a polyhedral chain $P$ and a $C^1$ diffeomorphism $f : \mathbf{R}^N \to \mathbf{R}^N$ such that $\mathrm{M}(A - f_\# P) < \epsilon$. (See [Fe, 4.1.28], or [Mo, 7.1] for a similar result.)

1.3. *Approximate derivatives.* Let $f : S \to \mathbf{R}^\ell$ be a function defined on a Lebesgue measurable subset $S$ of $\mathbf{R}^k$. We need the following facts about approximate differentiation [F, 3.1.4, 3.1.8]:

(1) $f$ is approximately differentiable almost everywhere if and only if there exist measurable sets $S_1 \subset S_2 \subset \ldots$ of $S$ such that $\mathcal{L}^k(S \setminus \cup S_i) = 0$ and such that the restriction $f|S_i$ of $f$ to $S_i$ is Lipschitz.

(2) $f$ is approximately differentiable almost everywhere if and only if for each $i = 1, \ldots, k$ and for almost every $x \in \mathbf{R}^k$ with $x_i = 0$, the function

$$t \mapsto f(x + t\mathbf{e}_i)$$

is approximately differentiable almost everywhere in $\{t : (x + t\mathbf{e}_i) \in S\}$ (where $\mathbf{e}_i$ is the $i^\mathrm{th}$ basis vector of $\mathbf{R}^k$).

For the purposes of this paper, one could take (1) as the definition of "approximately differentiable almost everywhere" and prove (2) as a theorem.

## 2. Flat chains of dimension $0$

2.1. THEOREM. *There is a homomorphism* $\chi : \mathcal{F}_0(\mathbf{R}^N; \mathbf{G}) \to \mathbf{G}$ *such that*

(1) $\chi(\sum g_i[x_i]) = \sum g_i$.
(2) $\chi(\partial T) = 0$ *for every* $1$-*chain* $T$.
(3) $|\chi(A)| \leq \mathcal{F}(A)$.
(4) $\mathcal{F}(A) \leq |\chi(A)| + \mathrm{M}(A)\,\mathrm{diam}(\mathrm{spt}\,A)$.

Of course $\mathcal{F}(A) \leq \mathrm{M}(A)$ for any flat chain $A$, so in particular (3) implies that $|\chi(A)| \leq \mathrm{M}(A)$.

*Remark.* The boundary of a 0-dimensional chain $A$ is by definition 0. However, in some ways it would make more sense to define the boundary of $A$ to be $\chi(A)$.

*Proof.* Define $\chi$ on polyhedral 0-chains by (1). Then (2) is obvious for polyhedral 1-chains $T$. Note by (1) that $|\chi(A)| \leq \mathrm{M}(A)$ for polyhedral 0-chains $A$. Thus if $T$ is a polyhedral 1-chain and $B = A - \partial T$, then

$$|\chi(A)| = |\chi(B) + \chi(\partial T)| = |\chi(B)| \leq \mathrm{M}(B).$$

Taking the infimum over all such $T$, we see that (3) holds for polyhedral 0-chains. This implies that $\chi$ has a unique $\mathcal{F}$-continuous extension to $\mathcal{F}_0(\mathbf{R}^N; \mathbf{G})$ since $\mathcal{P}_0$ is $\mathcal{F}$-dense in $\mathcal{F}_0$. Assertions (2) and (3) follow by continuity from the corresponding facts for polyhedral chains.

If $A = \sum g_i[x_i]$, let $x$ be a point in the convex hull of $\mathrm{spt}\,A$ and let $C$ be the cone over $A$ with vertex at $x$:

$$C = \sum g_i[x, x_i].$$

(Here $[x, x_i]$ is the oriented segment from $x$ to $x_i$.) Then $\partial C = A - \chi(A)[x]$, so

$$\mathcal{F}(A) \leq \mathrm{M}(\chi(A)[x]) + \mathrm{M}(C) \leq |\chi(A)| + \mathrm{M}(A)\,\mathrm{diam}(\mathrm{spt}\,A)$$

which proves (4) for polyhedral chains $A$. By continuity, it must hold for all 0-chains $A$. □

Let $\nu$ be a $\mathbf{G}$-valued Borel measure on $\mathbf{R}^N$. Recall that the total variation measure $|\nu|$ is given by

(*) $$|\nu|(S) = \sup \sum_i |\nu S_i|$$

where the supremum is over all partitions $S = \cup S_i$ of $S$ into countably-many Borel sets. The total variation of $\nu$ is $|\nu|(\mathbf{R}^N)$.

2.2. THEOREM. *Let $\mathcal{M}_0 = \mathcal{M}_0(\mathbf{R}^N; \mathbf{G})$ be the normed group of finite-mass flat 0-chains in $\mathbf{R}^N$ with coefficients in $\mathbf{G}$. There is an isomorphism $A \mapsto \psi_A$ from $\mathcal{M}_0$ to the group of compactly-supported $\mathbf{G}$-valued Borel measures on $\mathbf{R}^N$ of finite total variation (the norm of a measure being its total variation). Furthermore, $\mu_A = |\psi_A|$, where $|\psi_A|$ is the total variation measure associated to $\psi_A$.*

*Proof.* For Borel sets $S$, let

$$\psi_A(S) = \chi(A \llcorner S).$$

That $\psi_A$ is a compactly-supported Borel measure of finite total variation follows immediately from Theorem 2.1 (since $S \mapsto A \llcorner S$ defines a compactly-supported $\mathcal{M}_0$-valued Borel measure of finite total variation [FL, §4]).

By (*) and Theorem 2.1,

$$|\psi_A|(S) = \sup \sum_i |\psi_A(S_i)|$$
$$= \sup \sum_i |\chi(A \llcorner S_i)|$$
$$\leq \sup \sum_i \mathrm{M}(A \llcorner S_i)$$
$$= \mu_A(S).$$

Thus $|\psi_A| \leq \mu_A$.

To see that $\psi$ is one-to-one, suppose $\psi_A = 0$. Partition $\mathbf{R}^N$ into disjoint sets $S_i$, each of diameter $\leq \epsilon$. Then by Theorem 2.1(4):

$$\mathcal{F}(A) \leq \sum \mathcal{F}(A \llcorner S_i)$$
$$\leq \sum \left(|\chi(A \llcorner S_i)| + \mathrm{M}(A \llcorner S_i)\epsilon\right)$$
$$= \sum |\psi_A(S_i)| + \mathrm{M}(A)\epsilon$$
$$= 0 + \mathrm{M}(A)\epsilon.$$

Since $\epsilon$ is arbitrary, $\mathcal{F}(A) = 0$ and so $A = 0$. Thus $\psi$ is $1-1$.

To see that $\psi$ is surjective, let $\nu$ be a compactly supported $\mathbf{G}$-valued Borel measure on $\mathbf{R}^N$ with finite total variation. For $n = 0, 1, 2, \ldots$, let $\mathcal{Q}(n)$ be the standard partition of $\mathbf{R}^N$ into disjoint "half-open" cubes of side length $2^{-n}$. In particular, $(0, 2^{-n}]^N$ is one of the cubes of $\mathcal{Q}(n)$. Let

$$A_n = \sum_{Q \in \mathcal{Q}(n)} \nu(Q)[z_Q]$$

where $z_Q$ is the center of $Q$.

Now if $Q$ is a cube in $\mathcal{Q}(n)$ (with $n > 0$), then there is a unique cube $Q' \in \mathcal{Q}(n-1)$ containing $Q$; let $z'_Q = z_{Q'}$ be its center. Let $T_n$ be the polyhedral 1-chain given by:

$$T_n = \sum_{Q \in \mathcal{Q}(n)} \nu(Q)[z'_Q, z_Q].$$

Then

$$\partial T_n = A_n - A_{n-1},$$

so

$$\begin{aligned}\mathcal{F}(A_n - A_{n-1}) &\leq \mathrm{M}(T_n) \\ &= \sum_{Q \in \mathcal{Q}(n)} |\nu(Q)| \, |z'_Q - z_Q| \\ &= \sum_{Q \in \mathcal{Q}(n)} |\nu(Q)| 2^{-n-1} N^{1/2} \\ &\leq \left(2^{-n-1} N^{1/2}\right) |\nu|(\mathbf{R}^N).\end{aligned}$$

Thus the $A_n$ converge to a limit $A$. Note for any $Q \in \mathcal{Q}(n)$,

$$A_i \llcorner Q \to A \llcorner Q.$$

Also, for $i > n$, $\chi(A_i \llcorner Q) = \nu(Q)$. Thus

$$\psi_A(Q) = \chi(A \llcorner Q) = \nu(Q).$$

Since the various $Q$'s generate the $\sigma$-algebra of all Borel sets, $\psi_A = \nu$. Hence $\psi$ is surjective.

From the definition of $A_n$, we see immediately that $\mathrm{M}(A_n) \leq |\nu|(\mathbf{R}^N) = |\psi_A|(\mathbf{R}^N)$. Thus by lower-semicontinuity,

$$\mu_A(\mathbf{R}^N) = \mathrm{M}(A) \leq |\psi_A|(\mathbf{R}^N).$$

But we already showed that $|\psi_A| \leq \mu_A$, so this implies $|\psi_A| = \mu_A$. $\square$

2.3. *Rectifiable 0-chains.* A rectifiable 0-chain $A$ is a finite-mass chain for which there is a countable set $S$ with $A = A \llcorner S$. Equivalently, $A$ is rectifiable if and only if it is of the form

(*) $$\sum_{i=1}^{p} g_i[x_i]$$

where $\sum |g_i| < \infty$. (Here $p$ may be finite or $\infty$, depending on $A$.)

It is convenient to have a canonical representation of the form (*). We achieve this as follows. We may assume $A \neq 0$. Now for each $\epsilon > 0$, there are only finitely many points $x$ with $\mu_A\{x\} > \epsilon$. Hence there is a finite set $V$ of points $x$ for which $\mu_A\{x\}$ is a maximum. Let $x_1 \in V$ be the first point in $V$ according to the lexicographical ordering of $\mathbf{R}^N$. (In this ordering, $x < y$ means that for some $i$, the $i^{\text{th}}$ coordinate of $x$ is less than the $i^{\text{th}}$ coordinate of $y$, and that for every $j < i$, the $j^{\text{th}}$ coordinates of $x$ and of $y$ are equal.) We let $g_1$ be the multiplicity of $x_1$: $g_1 = \chi(A \llcorner \{x_1\})$. To get $x_2$ and $g_2$, we apply the same procedure to $A - g_1[x_1]$, and so on.

## 3. Slicing

Let $\sigma$ be a $k$-dimensional oriented polyhedron in $\mathbf{R}^N$ and let $P$ be an oriented $m$-dimensional affine subspace in $\mathbf{R}^N$. We say that $\sigma$ is transverse to $P$ if for each $0 \leq j \leq k$ and for each $j$-dimensional face $F$ of $\sigma$,
$$\dim F \cap P \leq j + m - N.$$
Transversality implies that either $\sigma$ and $P$ are disjoint, or $\sigma \cap P$ is a $j+m-N$ dimensional polyhedron containing interior points of $\sigma$.

We say that a polyhedral $k$-chain $A$ is transverse to $P$ if it can be expressed as a combination $A = \sum g_i[\sigma_i]$ of polyhedra $\sigma_i$ each of which is transverse to $P$. For such a chain $A$, we define the slice of $A$ by $P$ to be the $(N-m)$-chain
$$A \cap P = \sum g_i[\sigma_i \cap P].$$

Here we must adopt a convention by which the orientations on a pair of transverse subspaces $L$ and $M$ induce an orientation on $L \cap M$. The choice is somewhat arbitrary, but for this paper it is convenient to orient $L \cap M$ so that the following holds. If $w_1, \ldots, w_c$ is a correctly oriented basis for $L \cap M$, $u_1, \ldots, u_a, w_1, \ldots, w_c$ is a correctly oriented basis for $L$, and $v_1, \ldots, v_b, w_1, \ldots, w_c$ is a correctly oriented basis for $M$, then $u_1, \ldots, u_a, v_1, \ldots, v_b, w_1, \ldots, w_c$ is a correctly oriented basis for $\mathbf{R}^N$.

Clearly,
$$(A + B) \cap P = (A \cap P) + (B \cap P)$$
if $A$ and $B$ are transverse to $P$. Also,

(1) $$\partial(A \cap P) = \pm(\partial A) \cap P$$

with the sign depending only $N$ and the dimensions of $A$ and $P$. (With the orientation convention suggested above, the sign is always $+$.)

Now suppose $L$ is an oriented subspace of $\mathbf{R}^N$ and that $\Pi : \mathbf{R}^N \to L$ is the orthogonal projection. For each $x \in L$, we orient $\Pi^{-1}x$ so that $L \oplus \Pi^{-1}x$ induces the standard orientation on $\mathbf{R}^N$.

If $A$ is a polyhedral chain, then $A$ is transverse to $\Pi^{-1}x$ for almost every $x \in L$, and
$$\int_L \mathrm{M}(A \cap \Pi^{-1}x)\, dx \leq \mathrm{M}(A)$$

from which it follows (by (1)) that

$$\text{(2)} \quad \int_L \mathcal{F}(A \cap \Pi^{-1}x) \, dx \leq \mathcal{F}(A).$$

Now for each flat chain $A$, choose a sequence of polyhedral chains $A_i$ that converge rapidly to $A$ in that $\sum \mathcal{F}(A_i - A) < \infty$, and so that $\mathrm{M}(A_i) \to \mathrm{M}(A)$. Note we can do this in such a way that $(\tau_y A)_i = \tau_y(A_i)$ for all $y \in \mathbf{R}^N$. If $A$ is polyhedral, we let $A_i \equiv A$.

Note that $\sum \mathcal{F}(A_i - A_{i+1}) < \infty$. Thus by (2), for almost every $x$,

$$\sum \mathcal{F}(A_i \cap \Pi^{-1}x - A_{i+1} \cap \Pi^{-1}x) < \infty$$

and hence $A_i \cap \Pi^{-1}x$ converges rapidly to a limit chain, which we define to be $A \cap \Pi^{-1}x$.

Note that if $A'_i$ is another sequence of polyhedral chains converging rapidly to $A$, then $A_1, A'_1, A_2, A'_2, \dots$ also converges rapidly to $A$, which implies that $A'_i \cap \Pi^{-1}x \to A \cap \Pi^{-1}x$ for almost every $x \in L$.

3.1. PROPOSITION.

(1) $\int_L \mathrm{M}(A \cap \Pi^{-1}x) \, dx \leq \mathrm{M}(A)$.
(2) $\int_L \mathcal{F}(A \cap \Pi^{-1}x) \, dx \leq \mathcal{F}(A)$.
(3) If $A$ has finite mass and if $S$ is a Borel set in $\mathbf{R}^N$, then

$$(A \llcorner S) \cap \Pi^{-1}x = (A \cap \Pi^{-1}x) \llcorner S$$

for almost every $x \in L$.

(4) If $L$ and $M$ are orthogonal oriented subspaces of $\mathbf{R}^N$, then

$$A \cap \Pi^{-1}_{L \oplus M}(x,y) = \pm (A \cap \Pi^{-1}_L x) \cap \Pi^{-1}_M y$$

for almost every $(x,y) \in L \oplus M$ (the sign depending only on $N$ and the dimensions of $A$, $L$, and $M$).

*Proof.* The first two assertions follow immediately via Fatou's lemma from the corresponding facts for polyhedral chains.

Now (3) holds trivially for any set polyhedral set $S$ such that $A_i \llcorner S$ converges rapidly to $A \llcorner S$. By (1), it also holds for $\mu_A$-limits of such sets. But by [FL, §4], every Borel set is such a limit.

Finally, (4) follows immediately from the definition. □

3.2. THEOREM. *Let $A$ be a flat $k$-chain in $\mathbf{R}^N$ such that for each coordinate $k$-plane $L$ and for almost every $x \in L$,*

$$A \cap \Pi^{-1}x = 0.$$

*Then $A = 0$.*

*Proof.* Let $A_i \to A$ rapidly (i.e., let $\sum \mathcal{F}(A_i - A) < \infty$). Let $Q$ be a $k$-cube with sides parallel to $k$ of the coordinate axes, and let $Q^*$ be the $(N-k)$-cube with sides parallel to the other $(N-k)$ coordinate axes such that $Q$ and $Q^*$ have the same center and the same side length. Let $P = P_Q$ be the $(N-k)$-plane containing $Q^*$.

By definition of slicing,
$$\tau_x A_i \cap P \to \tau_x A \cap P = 0$$
rapidly for almost every $x$ parallel to $Q$.

Hence
$$\tau_y \left(\tau_x A_i \cap P\right) \llcorner Q^* \to 0$$
for almost every $y$ parallel to $Q^*$ (by [FL, 2.1]), which implies that for almost every $z \in \mathbf{R}^N$,
$$(\tau_z A_i \cap P) \llcorner Q^* \to 0$$
and hence that $g(z, Q, i) \to 0$, where
$$g(z, Q, i) = \chi\left((\tau_z A_i \cap P) \llcorner Q^*\right).$$

But the approximating chain $P^\epsilon(\tau_z A_i)$ in the deformation theorem [W3, 1.1, 2.2] is given by
$$P^\epsilon(\tau_z A_i) = \sum_Q g(z, Q, i)[Q]$$
where the sum is over all $k$-cubes $Q$ in the standard size $\epsilon$-cubeulation of $\mathbf{R}^N$. We have shown
$$P^\epsilon(\tau_z A_i) \to 0 \quad \text{as } i \to \infty$$
for almost every $z \in \mathbf{R}^N$. It follows that $P^\epsilon(\tau_z A) = 0$ for almost every $z$, and hence ([W3, 1.2(4) or 1.3]) that $A = 0$. □

3.3. THEOREM. *Let $A$ be a finite-mass flat $k$-chain. Let $S$ be a countable union of Borel sets $S_i$ of finite $k$-dimensional Hausdorff measure and suppose $\mu_A(\mathbf{R}^N \setminus S) = 0$. Let $\Pi : \mathbf{R}^N \to L$ be the orthogonal projection onto a $k$-plane $L$. Then*
$$A \cap \Pi^{-1} x$$
*is a rectifiable 0-chain for almost every $x \in L$.*

*If the $S_i$'s have finite integral geometric measure, then the conclusion holds for almost every $k$-plane $L$.*

In particular, the hypothesis holds if $A$ is rectifiable.

*Proof.* Note
$$(*) \qquad \mathcal{H}^k(S_i) \geq \int_{x \in L} \mathcal{H}^0(S_i \cap \Pi^{-1} x) \, d\mathcal{H}^k x$$

(cf. [Ma, 7.7] or [F, 2.10.25]) so $S_i \cap \Pi^{-1}x$ is finite for almost every $x$. Thus $S \cap \Pi^{-1}x$ is countable for almost every $x$. But (for almost every $x$)

$$A \cap \Pi^{-1}x = (A \llcorner S) \cap \Pi^{-1}x = (A \cap \Pi^{-1}x) \llcorner S = (A \cap \Pi^{-1}x) \llcorner (S \cap \Pi^{-1}x).$$

Since $A \cap \Pi^{-1}x$ lives on the countable set $S \cap \Pi^{-1}x$, it is rectifiable.

If $S_i$ has finite integral geometric measure, then the right-hand side of (*) is finite for almost every $L$; for such $L$ the rest of the proof is as before. □

## 4. Reduction to the graph of a continuous function

PROPOSITION. *Let $A$ be a finite-mass flat $k$-chain in $\mathbf{R}^N$. Let $L$ be a $k$-plane and let $\Pi : \mathbf{R}^N \to L$ be the orthogonal projection. Let $S$ be the set of $x \in L$ such that the slice $A \cap \Pi^{-1}x$ is rectifiable and nonzero. If $S$ has positive measure, then there exist a compact subset $K$ of $L$, a measurable function $g : K \to \mathbf{G}$ and a continuous function $\Phi : K \to \mathbf{R}^N$ such that*

1. $\mathcal{H}^k(K) > 0$.
2. $\inf_{x \in K} |g(x)| > 0$.
3. $(A \llcorner \Phi(K)) \cap \Pi^{-1}x = g(x)[\Phi(x)]$ *for $\mathcal{H}^k$-almost every $x \in K$.*
4. *If $V \subset K$ is a Borel set with $\mathcal{H}^k V > 0$, then*

$$A \llcorner \Phi(V) \neq 0.$$

*Proof.* Define functions $g : S \to \mathbf{G}$ and $\Phi : L \to \mathbf{R}^N$ as follows. For each $x \in S$, let $g(x)[\Phi(x)]$ be the first term in the canonical representation (§2.3) of $A \cap \Pi^{-1}x$. Since $g(x)$ is never 0,

$$S = \cup S_n$$

where

$$S_n = \left\{ x \in S : |g(x)| \geq \frac{1}{n} \right\}.$$

Since $\mathcal{H}^k S > 0$, there is an $n$ such that $\mathcal{H}^k S_n > 0$.

By Lusin's theorem, there is a compact subset $K$ of $S_n$ such that $\Phi|K$ is continuous and $\mathcal{H}^k K > 0$.

This proves (1) and (2). To prove (3),

$$(A \llcorner \Phi(K)) \cap \Pi^{-1}x = (A \cap \Pi^{-1}x) \llcorner \Phi(K)$$

$$= g(x)[\Phi(x)]$$

by Theorem 3.1(3).

It remains only to prove (4). Let $V$ be a Borel subset of $K$. Then

$$\mathrm{M}(A \llcorner \Phi(V)) \geq \int_{x \in L} \mathrm{M}((A \llcorner \Phi(V)) \cap \Pi^{-1}x) \, d\mathcal{H}^k x$$

$$= \int_{x \in L} \mathrm{M}((A \cap \Pi^{-1}x) \llcorner \Phi(V)) \, d\mathcal{H}^k x$$

$$= \int_{x \in V} \mathrm{M}(g(x)[\Phi(x)]) \, d\mathcal{H}^k x$$

$$= \int_{x \in V} |g(x)| \, d\mathcal{H}^k x$$

$$\geq \frac{1}{n} \mathcal{H}^k V$$

which is positive if $\mathcal{H}^k V > 0$. This proves (4). □

## 5. The one-dimensional case

5.1. LEMMA. *Let $T$ be a one-dimensional flat chain of finite mass and finite boundary mass. Then*

$$\mu_T \mathbf{B}(a, R) \geq \int_{r=0}^{R} |\chi((\partial T) \llcorner \mathbf{B}(a, r))| \, dr$$

(*where $\chi : \mathcal{F}_0(\mathbf{R}^N; \mathbf{G}) \to \mathbf{G}$ is the homomorphism of §2*).

COROLLARY. *If $T$ is a 1-chain with $\partial T = g[a] + E$ and if $T$ and $E$ have finite mass, then*

$$\mu_T \mathbf{B}(a, R) \geq R(|g| - \mu_E \mathbf{B}(a, R))$$

*for all $R$. (Here $\mathbf{B}(a, R)$ is the* open *ball of radius $R$ about $a$.)*

*Proof.* First suppose $T$ is polyhedral. Then

$$\mu_T \mathbf{B}(a, R) \geq \int_0^R \mathrm{M}(T \cap \partial \mathbf{B}(a, r)) \, dr$$

$$= \int_0^R \mathrm{M}(\partial(T \llcorner \mathbf{B}(a, r)) - (\partial T) \llcorner \mathbf{B}(a, r)) \, dr$$

$$\geq \int_0^R |\chi(\partial(T \llcorner \mathbf{B}(a, r)) - (\partial T) \llcorner \mathbf{B}(a, r))| \, dr$$

$$= \int_0^R |0 - \chi((\partial T) \llcorner \mathbf{B}(a, r))| \, dr$$

which proves the polyhedral case.

For a general $T$, let $T_i$ be polyhedral chains with $\sum \mathcal{F}(T_i - T) < \infty$ and $\mathrm{M}(T_i) \to \mathrm{M}(T)$. Then $\mu_{T_i} \to \mu_T$ and $T_i \llcorner \mathbf{B}(a,r) \to T \llcorner \mathbf{B}(a,r)$ for almost all $r$ ([FL, §4.2]). Thus the inequality follows from the inequality for the $T_i$'s by Fatou's lemma. □

*Proof of corollary.*

$$\begin{aligned}
\mu_T \mathbf{B}(a, R) &\geq \int_0^R |\chi((\partial T) \llcorner \mathbf{B}(a,r))|\, dr \\
&\geq \int_0^R |\chi((g[a] + E) \llcorner \mathbf{B}(a,r))|\, dr \\
&\geq \int_0^R |g + \chi(E \llcorner \mathbf{B}(a,r))|\, dr \\
&\geq \int_0^R (|g| - |\chi(E \llcorner \mathbf{B}(a,r))|)\, dr \\
&\geq \int_0^R (|g| - |\mathrm{M}(E \llcorner \mathbf{B}(a,r))|)\, dr \\
&\geq R (|g| - \mathrm{M}(E \llcorner \mathbf{B}(a,R))).
\end{aligned}$$
□

5.2. PROPOSITION. *Let $A$ be a flat chain in $\mathbf{R}^N$. Then there is a rectifiable chain $C$ with $\partial C = -\partial A$. If $A$ is supported in a hyperplane, then we can choose $C$ so that $\mu_C(\mathrm{spt}\, A) = 0$.*

*Remark.* For many purposes, we may assume that $A$ is supported in a hyperplane, because the inclusion $\mathbf{R}^N \subset \mathbf{R}^{N+1}$ allows us to regard $A$ as a chain in $\mathbf{R}^{N+1}$.

*Proof.* Let $A_i$ be a sequence of polyhedral chains with $A_1 = 0$ and $\sum \mathcal{F}(A_i - A) < \infty$. If $A$ is supported in a hyperplane $H$, we may choose the $A_i$'s to lie in $H$. Let $v$ be a unit normal to $H$, and translate $A_i$ a small amount in the direction of $v$ to get a polyhedral chain $A_i'$ with $\mathcal{F}(A_i' - A) < 2\mathcal{F}(A_i - A)$. In particular, $\sum \mathcal{F}(A_i' - A_{i+1}') < \infty$. Thus there exist polyhedral $(k+1)$-chains $P_i$ with $\sum \mathrm{M}(A_i' - A_{i+1}' - \partial P_i) < \infty$. Note we can choose $P_i$ to be disjoint from $H$. Let $C = \sum (A_i' - A_{i+1}' - \partial P_i)$. Then $\partial C = \lim_{n \to \infty} \sum_{i=0}^n (\partial A_i' - \partial A_{i+1}') = \lim_{n \to \infty} (\partial A_0' - \partial A_{n+1}') = 0 - \partial A$. If $A$ does not lie in a hyperplane, we let $A_i' = A_i$ and argue as before. □

5.3. PROPOSITION. *Let $A$ be a flat 1-chain in $\mathbf{R}^N$. Let $L$ be a line in $\mathbf{R}^N$ and let $\Pi : \mathbf{R}^N \to L$ be the orthogonal projection.*

*Suppose there is a compact set $K \subset L$, a function $g : K \to \mathbf{G}$, and a continuous function $\Phi : K \to \mathbf{R}^N$ such that*

$$\inf_{x \in K} |g(x)| = \delta > 0$$

*and such that*

$$A \cap \Pi^{-1}x = g(x)[\Phi(x)]$$

*for almost every $x \in K$.*

*Then there exist measurable sets $V_1 \subset V_2 \subset \cdots \subset K$ such that $\mathcal{L}^1(K \setminus \cup V_i) = 0$ and such that $\Phi|V_i$ is Lipschitz for each $i$. In particular, $\Phi$ is approximately differentiable at almost every point of $K$.*

*Proof.* By Proposition 5.2, there is a rectifiable 1-chain $C$ with $\partial(A+C) = 0$ and $\mu_C(\text{spt } A) = 0$. Let $C(x) = C \cap \Pi^{-1}x$. Then for almost every $x \in K$,

$$\mu_{C(x)}\{\Phi(x)\} = 0.$$

Also, by the Radon-Nykodym theorem, for $\mathcal{L}^1$-almost every $x \in K$ we have

$$\frac{d\nu}{d\mathcal{L}^1}x < \infty$$

where $\nu$ is the measure on $L$ given by

$$\nu(S) = \mu_{A+C}(\Pi^{-1}S).$$

Thus for almost every $x \in K$, there is an $n < \infty$ such that

(1) $$\mu_{C(x)}\mathbf{B}(\Phi(x), \tfrac{3}{n}) < \frac{1}{3}\delta$$

and

(2) $$\frac{\mu_{A+C}\Pi^{-1}[x, x+t]}{t} \leq n \quad \text{for } 0 < t < \frac{1}{n}.$$

Let $K_n$ be the set of $x \in K$ for which (1) and (2) hold. Then $K_1 \subset K_2 \subset K_3 \subset \ldots$ and $\mathcal{L}^1(K \setminus \cup K_i) = 0$.

Now we claim that

$$\limsup_{x,y \in K_n, |x-y| \to 0} \frac{|\Phi(x) - \Phi(y)|}{|x-y|} < \infty.$$

For let $x, y \in K_n$ with $x < y$. By the continuity of $\Phi$, we know that for $|x-y|$ sufficiently small, we will have

(3) $$|\Phi(x) - \Phi(y)| < \frac{1}{n}.$$

In particular this implies that $|x - y| < \frac{1}{n}$.

For such $x$ and $y$, let $T = (A+C) \llcorner \Pi^{-1}[x, y]$. Then

$$\partial T = -(A+C) \cap \Pi^{-1}x + (A+C) \cap \Pi^{-1}y$$
$$= -g(x)[\Phi(x)] - C(x) + g(y)[\Phi(y)] + C(y).$$

Let $r = |\Phi(x) - \Phi(y)|$. By Corollary 5.1,

$$(4) \quad \mathrm{M}(T) \geq (|g(x)| - \mu_{C(x)}\mathbf{B}(\Phi(x),r) - \mu_{C(y)}\mathbf{B}(\Phi(x),r))r$$
$$\geq (|g(x)| - \mu_{C(x)}\mathbf{B}(\Phi(x),r) - \mu_{C(y)}\mathbf{B}(\Phi(y),2r))r$$
$$\geq (\delta - \delta/3 - \delta/3)r$$
$$= \frac{\delta}{3}|\Phi(x) - \Phi(y)|$$

by (1) and (3). On the other hand, by (2),

$$\mathrm{M}(T) \leq n|x-y|.$$

This together with (4) implies

$$|\Phi(x) - \Phi(y)| \leq \frac{3n}{\delta}|x-y|. \qquad \square$$

## 6. Proof of the rectifiable slices theorem

We have already proved (§3.3) that if $A$ is rectifiable, then almost all of its slices are rectifiable. Thus suppose $A$ is a nonzero finite mass chain almost all of whose 0-dimensional slices (by planes parallel to coordinate planes) are rectifiable. We need to show that $A$ is rectifiable.

Note it suffices to show that, for any such $A$, there is a countably rectifiable set $S$ with $\mu_A S > 0$. For among such sets $S$, there is one that maximizes $\mu_A S$. If $A$ were not equal to $A \llcorner S$, then there would be a countably rectifiable set $S'$ for which $A' \llcorner S' \neq 0$, where $A' = A - A \llcorner S$. But then $S \cup S'$ would have greater $\mu_A$-measure than $S$, contradicting the choice of $S$.

By Theorem 3.2, there must be at least one coordinate $k$-plane $L$ for which the hypotheses of Section 4 hold. By rotating $A$, we may assume that it is the plane spanned by the first $k$ basis vectors of $\mathbf{R}^N$. Let $K \subset \mathbf{R}^k$, $g$, and $\Phi$ be as in Proposition 4. Of course

$$\Pi : \mathbf{R}^N \to \mathbf{R}^k$$
$$\Pi(x_1, \ldots, x_N) = (x_1, \ldots, x_k).$$

If $k = 1$, then we are done by Proposition 5.3 (applied to $A \llcorner \Phi(K)$). Thus we assume $k > 1$. Let $\pi : \mathbf{R}^N \to \mathbf{R}$ and $\tilde{\pi} : \mathbf{R}^N \to \mathbf{R}^{k-1}$ be the projections given by

$$\pi(x_1, \ldots, x_N) = x_1$$
$$\tilde{\pi}(x_1, \ldots, x_N) = (x_2, \ldots, x_k).$$

For $y \in \mathbf{R}^{k-1}$, let $K_y = \{x : (x,y) \in K\}$ and let $A_y$ be the 1-chain given by

$$A_y = (A \llcorner \Phi(K)) \cap \tilde{\pi}^{-1} y.$$

Note that for almost every $(x, y) \in K \subset \mathbf{R} \times \mathbf{R}^{k-1}$,

$$
\begin{aligned}
(1) \qquad A_y \cap \pi^{-1}x &= \pm(A \llcorner \Phi(K)) \cap \Pi^{-1}(x, y) \\
&= g(x, y)[\Phi(x, y)].
\end{aligned}
$$

In particular, for almost every $y \in \mathbf{R}^{k-1}$, (1) will hold for almost every $x \in K_y$.

Thus by Proposition 5.3 applied to $A_y$, $\Phi(\cdot, y)$, and $g(\cdot, y)$, we see that $x \mapsto \Phi(x, y)$ is approximately differentiable for almost every $x \in K_y$.

Likewise, for each $i = 2, \ldots, k$ and for almost every $x \in \mathbf{R}^k$ with $x_i = 0$, the function

$$t \mapsto \Phi(x + t\mathbf{e}_i)$$

is approximately differentiable almost everywhere. Hence (§1.3) $\Phi$ is approximately differentiable almost everywhere in $K$, which implies that there is a Borel set $V \subset K$ with $\mathcal{H}^k V > 0$ such that the restriction of $\Phi$ to $V$ is Lipschitz.

But then $\Phi(V)$ is rectifiable and $\mu_A \Phi(V) > 0$ (by §4).

This completes the proof of the rectifiable slices theorem.

6.1. COROLLARY. *Let $A$ be a finite mass k-chain. If there is a Borel set $S$ of finite k-dimensional Hausdorff (or integral geometric) measure such that $A = A \llcorner S$, then $A$ is rectifiable.*

*Proof.* By Proposition 3.3, almost every 0-dimensional slice of $A$ is rectifiable. Hence $A$ is rectifiable. □

## 7. Groups for which every finite-mass flat chain is rectifiable

7.1. THEOREM. *Let $\mathbf{G}$ be a complete normed abelian group, and let $0 \leq k < N$. The following are equivalent:*

(1) *Every flat k-chain of finite mass in $\mathbf{R}^N$ is rectifiable.*
(2) *There is no nonconstant continuous path of finite length in $\mathbf{G}$.*

Of course the length $L(\gamma; a, b)$ of a continuous path $\gamma : [a, b] \to \mathbf{G}$ is the supremum of

$$\sum_{i=1}^{m} |\gamma(t_i) - \gamma(t_{i-1})|$$

over all $m$ and all partitions $a = t_0 \leq t_1 \leq \cdots \leq t_m = b$.

*Proof that* (1) $\implies$ (2). Suppose there is a nonconstant finite length path $\gamma : [a, b] \to \mathbf{G}$. By reparametrizing, we may assume that it is parametrized by arc-length:

$$L(\gamma; s, t) = t - s \qquad \text{for } a \leq s \leq t \leq b.$$

Let $\nu$ be the **G**-valued Borel measure supported in $[a,b]$ such that

$$\nu[s,t] = \gamma(t) - \gamma(s)$$

for $a \leq s \leq t \leq b$. Note that $|\nu| = \mathcal{H}^1 \llcorner [a,b]$.

Thus by subsection 2.2 there is a finite-mass 0-chain $A$ with $\psi_A = \nu$. Note that $A$ is not 0 and $A$ is not rectifiable since

$$\mu_A(\{x\}) = |\nu|\{x\} = 0$$

for every point $x$.

If $k > 0$, let $Q$ be an oriented $(k-1)$-dimensional cube in $\mathbf{R}^{N-1}$. Then $B = A \times Q$ is a finite-mass $k$-chain that is not rectifiable, since $\mu_B = \mathcal{H}^{k+1}\llcorner([a,b] \times Q)$. □

*Proof that* (2) $\implies$ (1). Suppose there exist nonrectifiable chains of finite mass, and let $N$ be the smallest integer such that $\mathbf{R}^N$ contains such a chain $A$. Then $A$ must be a 0-dimensional chain. For if it were a $k$-chain with $k > 0$, then by the rectifiable slices theorem (together with Proposition 3.1(1)), we could slice it by a suitable $(N-k)$ plane $P$ to get a nonrectifiable finite mass 0-chain in $P \cong \mathbf{R}^{N-k}$, contradicting the choice of $N$.

We may assume that $A$ is purely unrectifiable, i.e., that there are no points of positive $\mu_A$ measure. For if it is not, we can replace $A$ by

$$A - A \llcorner \{x : \mu_A(x) > 0\}.$$

It follows that $\mu_A H = 0$ for every hyperplane of $\mathbf{R}^n$, since otherwise $A \llcorner H$ would be a nonzero purely unrectifiable chain in $H \cong \mathbf{R}^{N-1}$, contradicting the choice of $N$.

Since $A$ is nonzero, there is a Borel set $S$ with $\psi_A(S) \neq 0$ (by §2.2). We may assume that $\psi_A(\mathbf{R}^N) \neq 0$, for if it is not, we can replace $A$ by $A \llcorner S$.

Let $L : \mathbf{R}^N \to \mathbf{R}$ be a nonconstant linear function and let

$$\gamma : [a,b] \to \mathbf{G}$$
$$\gamma(t) = \psi_A L^{-1}(-\infty, t]$$

where $[a,b]$ is an interval whose interior contains $L(\operatorname{spt} A)$.

Now $\gamma$ is nonconstant since $\gamma(a) = 0$ and $\gamma(b) = \psi_A(\operatorname{spt} A) = \psi_A \mathbf{R}^N \neq 0$.

For $s < t$,

$$|\gamma(s) - \gamma(t)| = |\psi_A(L^{-1}(s,t])|$$
$$\leq \mu_A L^{-1}(s,t].$$

This implies that

$$\lim_{s \to t^-} |\gamma(s) - \gamma(t)| \leq \mu_A L^{-1} s$$
$$\lim_{t \to s^+} |\gamma(s) - \gamma(t)| = 0.$$

But $\mu_A L^{-1} s = 0$ since $L^{-1} s$ is a hyperplane of $\mathbf{R}^N$. Thus $\gamma$ is a nonconstant continuous path. Finally, it has finite length since if $a = t_0 \leq \cdots \leq t_m = b$, then

$$\begin{aligned}
\sum |\gamma(t_i) - \gamma(t_{i-1})| &= \sum |\psi_A L^{-1}(t_{i-1}, t_i]| \\
&\leq \sum \mu_A L^{-1}(t_{i-1}, t_i] \\
&= \mu_A L^{-1}(a, b] \\
&= \mathrm{M}(A).
\end{aligned}$$ □

## 8. Size and other weighted area functionals

Mass is a weighted area, with densities given by the group norm of the coefficients. In some situations, it is useful to consider other weighted areas. Let $\phi : \mathbf{G} \to [0, \infty]$ be any function that is even (i.e., such that $\phi(g) \equiv \phi(-g)$), lower-semicontinuous, and sub-additive, and such that $\phi(0) = 0$. Then $\phi$ induces a lower-semicontinuous functional (which we also denote by $\phi$) on flat chains:

$$\phi(A) = \inf \{\liminf \phi(A_i) : A_i \to A, \ A_i \in \mathcal{P}_k(\mathbf{R}^N; \mathbf{G})\}$$

where $\phi$ is defined on polyhedral chains by

$$\phi : \sum g_i [\sigma_i] \mapsto \phi(g_i) \mathrm{area}(\sigma_i)$$

if the $\sigma_i$'s are nonoverlapping. See [W3, §6] for basic properties of such functionals.

8.1. THEOREM. *Let $\mathbf{G}$ be a complete normed abelian group, $\phi : \mathbf{G} \to [0, \infty]$ be as above, and $0 \leq k < N$. The following are equivalent*:

(1) *Every flat $k$-chain $A$ with $\phi(A) + \mathrm{M}(A) < \infty$ is rectifiable.*
(2) *There is no nonconstant path in $\mathbf{G}$ that is both continuous and of finite length with respect to the metric $d(g, h) = \phi(g - h) + |g - h|$.*

*Proof.* The proof is a straightforward modification of 7.1. □

One such $\phi : \mathbf{G} \to [0, \infty]$ of particular interest is

$$\phi_s(g) = \begin{cases} 1, & \text{if } g \neq 0 \\ 0, & \text{if } g = 0. \end{cases}$$

Let us call $\phi_s(A)$ the **flat size** of the chain $A$.

8.2. PROPOSITION. *Let $A$ be a flat $k$-chain of finite mass. Then the following three quantities are equal*:

(1) *The flat size of $A$,*
(2) *The infimum of $\mathcal{H}^k(S)$ among Borel sets $S$ for which $A = A \llcorner S$,*
(3) *The infimum of the $k$-dimensional integral geometric measure of Borel sets $S$ for which $A = A \llcorner S$.*

*If this common quantity is finite, then $A$ is rectifiable.*

The quantities (2) and (3) might be called Hausdorff size and integral geometric size, respectively. Note that they are defined only when $M(A) < \infty$. Since flat size is defined for all flat chains and agrees with the other two whenever they are defined, it makes sense simply to call it size.

*Proof.* Unless all three quantities are infinite, $A$ must be rectifiable (by 6.1 and 8.2). For rectifiable chains, the Hausdorff and integral geometric sizes clearly coincide, and are also equal to the flat size [W3, §6]. □

STANFORD UNIVERSITY, STANFORD CA
*E-mail address*: white@math.stanford.edu